\documentstyle[11pt]{article}

\begin{document}

\title{Preservation theorems for countable support forcing iterations}
\author{
Chaz Schlindwein \\
Division of Mathematics and Computer Science \\
Lander University \\
Greenwood, South Carolina 29649, USA\\
{\tt chaz@lander.edu}}

\maketitle

\def\implies{\Rightarrow}
\def\rk{{\rm rk}}
\def\overb{{\overline b}}
\def\overw{{\overline w}}
\def\overx{{\overline x}}
\def\overy{{\overline y}}
\def\overz{{\overline z}}
\def\hht{{\rm ht}}
\def\forces{\mathbin{\parallel\mkern-9mu-}}
\def\notforces{\,\nobreak\not\nobreak\!\nobreak\forces}
\def\Gen{{\rm Gen}}
\def\calm{{\cal M}}
\def\bfone{{\bf 1}}

\def\restr{\,\hbox{\vrule height8pt width.4pt depth0pt
\vrule height7.75pt width0.3pt depth-7.5pt\hskip-.2pt
\vrule height7.5pt width0.3pt depth-7.25pt\hskip-.2pt
\vrule height7.25pt width0.3pt depth-7pt\hskip-.2pt
\vrule height7pt width0.3pt depth-6.75pt\hskip-.2pt
\vrule height6.75pt width0.3pt depth-6.5pt\hskip-.2pt
\vrule height6.5pt width0.3pt depth-6.25pt\hskip-.2pt
\vrule height6.25pt width0.3pt depth-6pt\hskip-.2pt
\vrule height6pt width0.3pt depth-5.75pt\hskip-.2pt
\vrule height5.75pt width0.3pt depth-5.5pt\hskip-.2pt
\vrule height5.5pt width0.3pt depth-5.25pt}\,}

\def\overtau{{\overline\tau}}

\def\cf{{\rm cf}}
\def\sup{{\rm sup}}
\def\supt{{\rm supt}}
\def\dom{{\rm dom}}
\def\range{{\rm range}}
\def\calj{{\cal J}}
\def\calc{{\cal C}}

\section{Introduction}

A major theme of [12] is preservation theorems for iterated
forcing. These are theorems of the form ``if $\langle P_\xi\,\colon
\xi\leq\kappa\rangle$ is a countable support forcing iteration based on
$\langle\dot Q_\xi\,\colon\xi<\kappa\rangle$ and each $\dot Q_\xi$
has property such-and-such then $P_\kappa$
has property thus-and-so.'' The archetypal preservation theorem
is the Fundamental Theorem of Proper Forcing [12, chapter III],
which states that if each $\dot Q_\xi$ is proper in $V[G_{P_\xi}]$
then $P_\kappa$ is proper.
Typically, the property enjoyed by $P_\kappa$ ensures that
either $\omega_1$ is not collapsed, or that no new reals are
added. In this paper we introduce two preservation theorems, one for not
collapsing $\omega_1$ and one for not adding reals, which include
many of the preservation theorems of [12] as special cases.
We shall see that some preservation theorems from [12] which
use revised countable support iteration are true also for
countable support iteration; for example, we have that
any countable support iteration of semi-proper forcings is
 semi-proper. These results lessen the importance of
  the concept of revised countable
support iterations.  

For preserving $\omega_1$, we introduce the class of hemi-proper forcings.
This generalizes semi-properness, and, assuming CH, includes also
Namba forcing and the forcing $P[S]$ consisting of all 
increasing countable sequences from~$S$,
ordered by reverse end-extension,
 where $S\subseteq S^2_0=_{\rm def}
\{\alpha<\omega_2\,\colon\allowbreak
\cf(\alpha)=\omega\}$ and $S$ is stationary.
This is the forcing notion used by Shelah in his solution to
a problem of H. Friedman, namely the question of whether every
stationary subset of $S^2_0$ may contain a closed subset
of order-type $\omega_1$.

For not adding reals, we introduce $\mu$-pseudo-completeness
 (our terminology clashes with
the terminology of Shelah in this instance).
The hypothesis of our preservation theorem is that each
$\dot Q_\xi$ is  $\mu$-pseudo-complete relative to
$P_\xi$, and the conclusion is that $P_\kappa$ does not
add any elements of ${}^\omega\mu$. We show that under CH,
both Namba forcing and $P[S]$ satisfy the hypothesis of the
theorem; for Namba forcing we require $\mu\in\{2,\omega\}$ and
for $P[S]$ we allow any $\mu$. Also Prikry forcing
works for $\mu$ less than the measurable cardinal.

We also make an observation regarding the construction of models
in which the continuum is larger than $\aleph_2$ via
iterated forcing; this sort of construction, we observe,
is not so difficult as has been previously believed.

{\bf Notation.}  Our notation follows [10], except
as noted in definition 1. We set
$M^P$ equal to the set (or class) of $P$-names which are in $M$.
This is different from the class of names whose values are
forced to be in $M$, and it is different from the class of
names whose values are forced to be in $M[G_P]$ (of course, for
any $G_P$ and any name $\dot x$ which is forced to be in $M[G_P]$
there are $p\in G_P$ and $\dot y \in M^P$ such that $p\forces``\dot x
=\dot y$'').
The notation $\dot P_{\eta,\alpha}$ is used in the context
of a forcing iteration $\langle P_\beta\,\colon\beta\leq\alpha\rangle$
based on $\langle\dot Q_\beta\,\colon\beta<\alpha\rangle$; it
denotes a $P_\eta$-name characterized by
$V[G_{P_\eta}]\models``\dot P_{\eta,\alpha}=\{p\restr[\eta,\alpha)
\,\colon p\restr\eta\in G_{P_\eta}$ and $p\in P_\alpha\}$.''
By $p\restr[\eta,\alpha)$ we do not mean the check (with respect to
$P_\eta$) of the restriction of $p$ to the interval
$[\eta,\alpha)$, but rather we mean the $P_\eta$-name for the function
$f$ with domain equal to $[\check\eta,\check\alpha)$
such that $f(\check\beta)$ is the $P_\eta$-name for the
$\dot P_{\eta,\beta}$-name corresponding to the $P_\beta$-name
$p(\beta)$ (see  [12, definition II.2.3, page 45]).
We shall use such facts as $\bfone\forces``\dot P_{\eta,
\alpha}$ is a poset;'' see [10] for a proof.
We shall use the notation of [6] regarding forcing names; see
[6, page 188].

\section{Hemi-properness}

In this section we introduce hemi-properness and show that
hemi-properness is preserved under countable support iterations.
As a warm-up, we show that semi-properness is 
preserved under countable support
iterations.

\proclaim Definition 1. Suppose\/ $\langle P_\xi\,\colon\xi\leq
\alpha\rangle$ is a countable support forcing iteration
and\/ $\eta<\alpha$.
 By\/ $\dot P_{\eta,\alpha}^M$
we mean the name\/ $\{\langle p\restr[\eta,\alpha),p\restr\eta\rangle
\,\colon p\in P_\alpha$ {\rm and $p\restr[\eta,\alpha)\in M^{P_\eta}\}$.}
This is in contrast with the object\/ $\{\langle \dot q
,p\rangle\,\colon p\in P_\eta$ {\rm and
$p\forces``\dot q\in \dot P_{\eta,\alpha}$'' and
$\dot q\in M^{P_\eta}\}$.}
Notice that the assertion\/ {\rm $r\forces_{P_\eta}``\dot s\in\dot 
P_{\eta,\alpha}^M
$''} is stronger than the assertion\/
{\rm $r\forces``\dot s\in\dot P_{\eta,\alpha}\cap
M[G_{P_\eta}]$''} and it is also different from\/
{\rm $r\forces``\dot s\in\dot P_{\eta,\alpha}\cap \check M$.''}
By\/ {\rm ``$\supt(p)$''} we mean\/
$\{\beta\in\dom(p)\,\colon p(\beta)\ne\bfone_{\dot Q_\beta}\}$.
This is in contrast with\/ {\rm [10],} where\/ {\rm $
``\supt(p)$''} was used to mean\/ {\rm $\{\beta\in\dom(p)\,\colon
p\restr\beta\notforces``p(\beta)=\bfone_{\dot Q_\beta}$''$\}$.}

In the following three lemmas we establish the basic
  facts about $\dot P^M_{\eta,\alpha}$.

\proclaim Lemma 2. Suppose $\langle P_\xi\,\colon\xi\leq\kappa\rangle$ is
a countable support forcing iteration and\/ $M$ is a countable
elementary submodel of\/ $H_\lambda$ for some sufficiently
large regular\/ $\lambda$ and\/ $N$ is
a countable elementary substructure of\/ $H_\theta$,
where $\theta$ is a regular cardinal greater than\/
$\lambda$, and suppose $P_\kappa\in M\in N$.
Suppose $\sigma\in N^{P_\kappa}$ and    $\eta\in\kappa\cap M$ and\/
$p\in P_\eta$ and\/ $\dot q\in N^{P_\eta}$ and\/
 $p\forces``\sigma\in M[G_{P_\eta}]^{\dot P_{\eta,\kappa}}$ {\rm
 and $\dot q\in\dot P_{\eta,\kappa}^M$.''}
Then there is $\dot r\in N^{P_\eta}$ and\/ $\tau\in N^{P_\eta}$
such that\/ {\rm
$p\forces``\dot r\in\dot P^M_{\eta,\kappa}$ and $\dot r\leq\dot q$
and $\dot r\forces_{\dot P_{\eta,\kappa}}`$if $\sigma$ is
an ordinal then $\sigma=\check\tau$' and $\tau\in
M[G_{P_\eta}]$.''}

Proof: Work in $N$ (so we must not
refer to $p$, which is not in $N$).
Take $\calj\subseteq\{s\in P_\eta\,\colon(\exists p(s)\in P_\kappa)
(\exists\sigma(s)\in M^{P_\eta})
(p(s)\restr\eta=s$ and $p(s)\restr[\eta,\kappa)\in M^{P_\eta}$
and $s\forces``$if $\dot q\in
\dot P^M_{\eta,\kappa}$ and $\sigma\in
M[G_{P_\eta}]^{\dot P_{\eta,\kappa}}$
then $p(s)\restr[\eta,\kappa)=\dot q$ and $\sigma=\sigma(s)$''$
)\}$ a maximal antichain.
For each $s\in\calj$ take $\dot r(s)\in M^{P_\eta}$ and $\tau(s)\in
M^{P_\eta}$ such that $\bfone\forces_{P_\eta}``$if
$\dot q\in\dot P^M_{\eta,\kappa}$ and $\sigma\in M[G_{P_\eta}]^{\dot
P_{\eta,\kappa}}$ then $\dot r(s)\leq
p(s)\restr[\eta,\kappa)$ and $\dot r(s)\forces`$if $\sigma(s)$
is an ordinal then $\sigma(s)=\tau(s)$' and $\tau(s)$ is
an ordinal.'' For each $s\in\calj$ take $\calj(s)\in M$
such that $\calj(s)$ is a maximal antichain of $\{s'\in P_\eta\,\colon
(\exists p'(s')\in P_\kappa)(p'(s')\restr\eta=s'$
and $s'\forces``p'(s')\restr[\eta,\kappa)=\dot r(s)$''$)\}$.
Take $p^*(s)\in M^{P_\eta}$ to be forced to be a function with
domain equal to the interval $[\check\eta,\check\kappa)$
such that $(\forall\gamma\in[\eta,\kappa))(\forall s'\in
\calj(s))(s'\forces``p^*(s)(\check\gamma)=p'(s')(\check\gamma)$'').
Now take $\dot r=\{\langle p^*(s),s\rangle\,\colon s\in\calj\}$
and take $\tau\in V^{P_\eta}$ such that $(\forall s\in\calj)
(s\forces``\tau=\tau(s)$''). We are done.

\proclaim Lemma 3. Suppose\/
 $\langle P_\xi\,\colon\xi\leq\kappa\rangle$ is
a countable support forcing iteration  and\/ $\lambda$ is a
sufficiently large cardinal and\/ $M$ is a countable elementary
substructure of\/ $H_\lambda$. Suppose\/ $\eta\in \kappa\cap M$
and\/ $\{P_\kappa,\eta\}\subseteq M$. Suppose also\/
$p\in P_\eta$ and\/ {\rm $p\forces``\dot q\in\dot P_{\eta,\kappa}^M
$.''} Then\/ {\rm $p\forces``\supt(\dot q)
\subseteq \check M$.''}

Proof: Recall by definition 1 that we are using ``$\supt$''
in the sense of [6] rather than in the sense of [10].
Given $r\leq p$, take $s\in P_\kappa$ such that
$s\restr\eta\leq r$ and $s\restr\eta\forces``\dot q=
s\restr[\eta,\kappa)$'' and $s\restr[\eta,\kappa)\in M^{P_\eta}$. 
We have $s\restr\eta\forces``\supt(\dot q)\subseteq
\supt(s\restr[\eta,\kappa))\subseteq\check X$'' where
$X=\supt(s)$. Now, $X\cap[\eta,\kappa)$ is a countable
element of $M$, hence it is a subset of $M$. We are done.

\proclaim Corollary 4. Suppose\/ $\langle P_\alpha\,\colon\alpha\leq
\kappa\rangle$ is a countable support iteration and\/ $M$ is
a countable elementary submodel of\/ $H_\lambda$ and\/
$\eta<\kappa$ and\/ $\{P_\kappa,\eta\}\subseteq M$. Suppose\/
$p\in P_\eta$ and\/ {\rm $p\forces``\dot q\in\dot P_{\eta,\kappa}^M$.''}
Then there is\/ $r\in P_\kappa$ such that\/ $r\restr\eta=p$ and\/ {\rm
$p\forces``r\restr[\eta,\kappa)=\dot q$.''}

Remark: Notice that this is false if we weaken the
hypothesis $``p\forces`\dot q\in \dot P^M_{\eta,\kappa}$'\thinspace''
to $``p\forces`\dot q\in\dot P_{\eta,\alpha}\cap 
M[G_{P_\eta}]$.'\thinspace''

Observation: It has been remarked by Roitman that there is a
difficulty in obtaining a model in which the continuum is
larger than $\aleph_2$ using countable support iterations
(see [1, page 56]).  Notice that this difficulty disappears
with our approach.  The difficulty arises because CH fails in some
intermediate model, and the traditional method of proving
preservation theorems for countable support involves examining
the behavior of $\dot P_{\eta,\kappa}$ in the model
$V[G_{P_\eta}]$.  However, in our approach, so long as CH holds
in the ground model, the fact that it fails in intermediate
models is immaterial.

\proclaim Lemma 5. Suppose\/ $\eta<\beta<\alpha$ and\/ {\rm
$p\forces_{P_\eta}``\dot q\in\dot P^M_{\eta,\alpha}$''}
and\/ $\{\eta,\beta,P_\alpha\}\subseteq M$, where $M$ is
a countable elementary substructure of\/ $H_\lambda$.
Then\/ {\rm
 $p\forces_{P_\eta}``\dot q\restr\beta\forces_{\dot P_{\eta,\beta}}
`\dot q\restr[\beta,\alpha)\in\dot P^M_{\beta,\alpha}$.'\thinspace''
}

Proof: Given $p'\leq p$ take $r\in P_\alpha$
such that $r\restr\eta\leq p'$ and $r\restr[\eta,\alpha)\in M^{P_\eta}$
and $r\restr\eta\forces``r\restr[\eta,\alpha)=\dot q$.''
Then we have $r\restr\beta\forces``\dot q\restr[\beta,\alpha)
=r\restr[\beta,\alpha)\in \dot P^M_{\beta,\alpha}$.''
Hence $r\restr\eta\forces``r\restr[\eta,\beta)\forces`\dot q\restr[\beta,
\alpha)\in\dot P_{\beta,\alpha}^M$'\thinspace'' by [10, section 3].
Because $p'$ was an arbitrary condition below $p$ we are done.

We assume familiarity with the definition
of ``semi-proper'' from [12, chapter X].
Here we introduce the appropriate induction hypothesis
 for showing that semi-properness
 is preserved under countable support iterations.

\proclaim Definition 6. Suppose\/ $\langle P_\xi\,\colon\xi\leq\alpha
\rangle $ is a countable support iteration. We say that $P_\alpha$ is\/
{\rm strictly semi-proper} iff
 whenever $\lambda$ is a sufficiently large
regular
cardinal and $\eta<\alpha$ and $M$ is a countable elementary
substructure of $H_{\lambda}$ 
and $\{P_\alpha,\eta\}\subseteq M$ and
$\dot q\in V^{P_\eta}$ and\/
$p\in P_\eta$ and\/ $p\forces``\omega_1\cap M[G_{P_\eta}]=
\omega_1\cap M$ {\rm and
$\dot q\in\dot P_{\eta,\alpha}^M$,''} then there is\/
 $r\in P_\alpha$ such that\/
$r\restr\eta=p$ and\/ $p\forces``r\restr[\eta,\alpha)\leq\dot q$''
and $r\forces``\omega_1\cap M[G_{P_\alpha}]=\omega_1\cap M$.''

\proclaim Lemma 7. Suppose\/
 $\langle P_\xi\,\colon\xi\leq\kappa\rangle$ is a countable support
iteration
based on $\langle\dot Q_\xi\,\colon\xi<\kappa\rangle$ and
$P_\alpha$ is strictly semi-proper for every $\alpha<\kappa$,
and if $\kappa=\gamma+1$ then\/ 
{\rm $\bfone\forces_{P_\gamma}``\dot Q_\gamma$ is
semi-proper.''}
 Then
$P_\kappa$ is strictly semi-proper.

Proof: Suppose $\lambda$, $\eta$, $M$, $\dot q$, and
$p$ are as in definition 6. 

Case 1: $\gamma+1=\kappa$ for some $\gamma$.

Take $r_0\in P_\gamma$ such that
$r_0\restr\eta=p$ and $p\forces``r_0\restr[\eta,\gamma)\leq
\dot q\restr\gamma$'' and 
$r_0\forces``\omega_1\cap M[G_{P_\gamma}]=\omega_1\cap M$.''
Because $\bfone\forces_{P_\gamma}``\dot Q_\gamma$ is
semi-proper,'' we may take $\dot s\in V^{P_\gamma}$
such that $\bfone\forces``\dot s\leq\dot q(\gamma)$ and
$\dot s\forces``\omega_1\cap M[G_{P_\gamma}][G_{\dot Q_\gamma}]
=\omega_1\cap M[G_{P_\gamma}]$.''
Take $r\in P_\kappa$ such that $r\restr\gamma=r_0$ and
$r(\gamma)=\dot s$. Then
$r\forces``\omega_1\cap M[G_{P_\kappa}]=
\omega_1\cap M[G_{P_\gamma}]=\omega_1\cap M$'' and
we are done.

Case 2: $\kappa$ is a limit ordinal.

Let $\alpha=\sup(\kappa\cap M)$ and let $\langle \sigma_n
\,\colon n\in\omega\rangle$ enumerate $M^{P_\kappa}$.
Let $\langle \alpha_n\,\colon n\in\omega\rangle$ be an increasing
sequence of ordinals from $\alpha\cap M$ cofinal in $\alpha$
such that $\alpha_0=\eta$.

Build $\langle p_n,\dot q_n,\tau_n\,\colon n\in\omega\rangle$
such that $p_0=p$ and $\dot q_0=\dot q$ and
both of the following:

($i$) $p_n\in P_{\alpha_n}$ and $p_{n+1}\restr\alpha_n=p_n$
and $p_n\forces``p_{n+1}\restr[\alpha_n,\alpha_{n+1})\leq
\dot q_n\restr\alpha_{n+1}$'' and $p_n\forces``\omega_1\cap
M[G_{P_{\alpha_n}}]=\omega_1\cap M$''

($ii$) $p_{n+1}\forces``\dot q_{n+1}\in
\dot P^{M}_{\alpha_{n+1},\kappa}$ and
$\dot q_{n+1}\leq\dot q_n\restr[\alpha_{n+1},\kappa)$ and
$\dot q_{n+1}\forces`$if $\sigma_n\in\omega_1$
then $\sigma_n=\check\tau_n$' and $\tau_n\in
\omega_1\cap M[G_{P_{\alpha_{n+1}}}]$.''

This is possible by lemmas 2 and 5 and
the induction hypothesis.

Now take $r\in P_\kappa$
such that $\supt(r)\subseteq\alpha$ and
for all $n\in\omega$ we have $r\restr\alpha_n=p_n$.
By lemma 3 we have $p_{n+1}\forces``
\supt(\dot q_{n+1})\subseteq\alpha$'' and
hence we have $r\restr\alpha_{n+1}\forces``r\restr[\alpha_{n+1}
,\kappa)\leq\dot q_{n+1}$.'' Therefore
we have $r\restr\alpha_{n+1}\forces``r\restr[\alpha_{n+1},
\kappa)\forces`$if $\sigma_n\in\omega_1$ then
$\sigma_n\in \omega_1\cap M[G_{P_{\alpha_{n+1}}}]
=\omega_1\cap M$.'\thinspace''
Now suppose $r_1\leq $ and $r_1\forces``\sigma\in \omega_1\cap
M[G_{P_\kappa}]$. Take $r_2\leq r_1$ and $n\in\omega$ such that
$r_2\forces``\sigma=\sigma_n$.'' Then $r_2\forces``\sigma\in\check M$.''
We conclude
that $r\forces``\omega_1\cap M[G_{P_\kappa}]=
\omega_1\cap M$.''
The lemma is established.

\proclaim Lemma 8. Suppose\/ $P_\kappa$ is strictly
semi-proper. Then\/ $P_\kappa$ is semi-proper.

Proof: Given $\lambda$ sufficiently large, regular, and
$M$ a countable elementary submodel
of $H_{\lambda^{+}}$, with $P_\kappa\in M$
and $q\in P_\kappa\cap M$, take
$\eta=0$ in definition~6.
We obtain by lemma 7 a condition $r\leq q$
such that $r\forces``\omega_1\cap M[G_{P_\kappa}]
=\omega_1\cap M$.''
We are done.

\proclaim Theorem 9. Suppose\/ $\langle P_\xi\,\colon\xi\leq\kappa\rangle$
is a countable support iteration based on\/ $\langle\dot Q_\xi\,\colon
\xi<\kappa\rangle$, and for every\/ $\xi<\kappa$ we have
that\/ {\rm $\bfone\forces_{P_\xi}``
\dot Q_\xi$ is semi-proper.''}
Then\/ $P_\kappa$ is semi-proper.

Proof. By lemmas 7 and 8.

We now turn our attention to hemi-properness.

\proclaim Definition 10. We say that a poset\/ $P$ is {\rm hemi-proper}
iff whenever\/ $\lambda$ is an appropriately large regular cardinal
and\/ $M$ and\/ $N$ are countable elementary substructures of\/ $H_\lambda$
and\/ $P\in M\in N$ and $q\in P\cap M$, then\/  {\rm
$q\notforces``\omega_1\cap M[G_P]>\omega_1\cap N$.''}

\proclaim Lemma 11. Suppose\/ $P$ is hemi-proper. Then\/ $P$ does
not collapse\/ $\omega_1$.

Proof: Suppose $q\forces_{P}``\dot f$ maps $\omega$ onto $\omega_1^V
$.'' Take $M$ and $N$ countable elementary substructures of some
appropriate $H_{\lambda}$ such that $\{P,q,\dot f\}\subseteq M\in N$.
We have $q\notforces``\omega_1\cap M[G_{P}]>\omega_1\cap N$.''
Hence $q\notforces``\omega_1\subseteq M[G_{P}]$,'' a contradiction.

\proclaim Lemma 12. If\/ $P$ is semi-proper, then\/ $P$ is hemi-proper.

Proof: Immediate.

In theorem 19 we show that the class of hemi-proper forcings also
contains Namba forcing, assuming CH, and in theorem 22 we show that
$P[S]$ is hemi-proper for $S\subseteq S^2_0$ stationary.

It is instructive to see why the poset which adds a closed
unbounded subset to a given stationary subset $S$
of $\omega_1$
using countable conditions is not hemi-proper (if it were,
then our preservation theorem
would contradict the fact that one can collapse
$\omega_1$ with an $\omega$-length forcing iteration
built of posets of this form).  The reason is that $M$ and
$N$ may be such that $S\cap N\cap\{\omega_1\cap M'\,\colon
M'\prec H_\lambda\}\subseteq M$.

\proclaim Definition 13. Suppose\/ $\langle P_\xi\,\colon\xi\leq\alpha
\rangle $ is a countable support iteration. We say that\/ $P_\alpha$ is\/
{\rm strictly hemi-proper} iff whenever\/
 $\lambda$ is a sufficiently large
regular
cardinal and\/ $\eta<\alpha$ and\/ $M$ is a countable elementary
substructure of\/ $H_{\lambda^{+\eta+1}}$ and\/ $N$ is a countable
elementary substructure of\/ $H_{\lambda^{+\alpha+1}}$
and\/ $\{P_\alpha,\eta\}\subseteq M\in N$ and\/
$\dot q\in N^{P_\eta}$ 
and\/ $p\in P_\eta$ and\/ $p\forces``\omega_1\cap M[G_{P_\eta}]\leq
\omega_1\cap N$ {\rm and
$\dot q\in\dot P^M_{\eta,\alpha}$,''} then there is\/
 $r\in P_\alpha$ such that\/
$r\restr\eta=p$ and\/ $p\forces``r\restr[\eta,\alpha)\leq\dot q$''
and\/ {\rm $r\forces``\omega_1\cap M[G_{P_\alpha}]\leq\omega_1\cap N$.''}

\proclaim Lemma 14. Suppose\/
$\langle P_\xi\,\colon\xi\leq\kappa\rangle$ is a countable support
iteration
based on\/ $\langle\dot Q_\xi\,\colon\xi<\kappa\rangle$ and\/
$P_\alpha$ is strictly hemi-proper for every\/ $\alpha<\kappa$,
and if\/ $\kappa=\gamma+1$ then\/ 
{\rm $\bfone\forces_{P_\gamma}``\dot Q_\gamma$ is
hemi-proper.''}
 Then\/
$P_\kappa$ is strictly hemi-proper.

Proof: Suppose $\lambda$, $\eta$, $M$, $N$, $\dot q$, and
$p$ are as in definition 13. 

Case 1: $\gamma+1=\kappa$ for some $\gamma$.

Take $M'$ a countable elementary substructure of $H_{\lambda^{
+\kappa}}$ such that $M\in M'\in N$. Take
$r_0\in P_\gamma$ such that
$r_0\restr\eta=p$ and $p\forces``r_0\restr[\eta,\gamma)\leq
\dot q\restr\gamma$'' and $r_0\forces``\omega_1\cap M'[G_{P_\gamma}]\leq\omega_1\cap N$.''
Because $\bfone\forces_{P_\gamma}``\dot Q_\gamma$ is
hemi-proper,'' we may take $\dot s\in V^{P_\gamma}$
such that $\bfone\forces``\dot s\leq\dot q(\gamma)$ and
$\dot s\forces``\omega_1\cap M[G_{P_\gamma}][G_{\dot Q_\gamma}]
\leq\omega_1\cap M'[G_{P_\gamma}]$.''
Take $r\in P_\kappa$ such that $r\restr\gamma=r_0$ and
$r(\gamma)=\dot s$. Then
$r\forces``\omega_1\cap M[G_{P_\kappa}]\leq
\omega_1\cap M'[G_{P_\gamma}]\leq\omega_1\cap N$'' and
we are done.

Case 2: $\kappa$ is a limit ordinal.

Let $\alpha=\sup(\kappa\cap M)$ and let $\langle \sigma_n
\,\colon n\in\omega\rangle$ enumerate $\{\sigma\in N^{P_\kappa}\,\colon
\bfone\forces_{P_\kappa}``\sigma\in\omega_1\cap M[G_{P_\kappa}]$''$\}$.
Let $\langle \alpha_n\,\colon n\in\omega\rangle$ be an increasing
sequence of ordinals from $\alpha\cap M$ cofinal in $\alpha$
such that $\alpha_0=\eta$.
Let $M_n=N\cap H_{\lambda^{+\alpha_n+1}}$.

Build $\langle p_n,\dot q_n,\tau_n\,\colon n\in\omega\rangle$
such that $p_0=p$ and $\dot q_0=\dot q$ and
both of the following:

($i$) $p_n\in P_{\alpha_n}$ and $p_{n+1}\restr\alpha_n=p_n$
and $p_n\forces``p_{n+1}\restr[\alpha_n,\alpha_{n+1})\leq
\dot q_n\restr\alpha_{n+1}$'' and $p_n\forces``\omega_1\cap
M[G_{P_{\alpha_n}}]\leq\omega_1\cap M_n$''

($ii$) $p_{n+1}\forces``\dot q_{n+1}\in
\dot P^{M}_{\alpha_{n+1},\kappa}$ and
$\dot q_{n+1}\leq\dot q_n\restr[\alpha_{n+1},\kappa)$ and
$\dot q_{n+1}\forces`
\sigma_n=\check\tau_n$' and $\tau_n\in M[G_{P_{\alpha_{n+1}}}]$''
and $\dot q_{n+1}\in N^{P_{\alpha_{n+1}}}$

This is possible by lemmas 2 and 5 and
the induction hypothesis.

Now take $r\in P_\kappa$
such that $\supt(r)\subseteq\alpha$ and
for all $n\in\omega$ we have $r\restr\alpha_n=p_n$.
By lemma 3 we have $p_{n+1}\forces``
\supt(\dot q_{n+1})\subseteq\alpha$'' and
hence we have $r\restr\alpha_{n+1}\forces``r\restr[\alpha_{n+1}
,\kappa)\leq\dot q_{n+1}$.'' Therefore
we have $r\restr\alpha_{n+1}\forces``r\restr[\alpha_{n+1},
\kappa)\forces`
\sigma_n\in \omega_1\cap M[G_{P_{\alpha_{n+1}}}]
\subseteq \omega_1\cap M_{n+1}=\omega_1\cap N$.'\thinspace''
We conclude
that $r\forces``\omega_1\cap M[G_{P_\kappa}]\subseteq
\omega_1\cap N$.''
The lemma is established.

\proclaim Lemma 15. Suppose\/ $P_\kappa$ is strictly
hemi-proper. Then\/ $ P_\kappa$ is hemi-proper.

Proof: Given $\lambda$ sufficiently large, regular, and
$M'$ and $N$ countable elementary submodels
of $H_{\lambda^{+\kappa+1}}$, with $P_\kappa\in M'\in N$
and $q\in P_\kappa\cap M'$, take
$\eta=0$ and $M=M'\cap H_{\lambda^+}$ in definition 13.
We obtain by lemma 14 a condition $r\leq q$
such that $r\forces``\omega_1\cap M'[G_{P_\kappa}]
=\omega_1\cap M[G_{P_\kappa}]\leq\omega_1\cap N$.''
We are done.

\proclaim Theorem 16. Suppose\/ $\langle P_\xi\,\colon\xi\leq\kappa\rangle$
is a countable support iteration based on\/ $\langle\dot Q_\xi\,\colon
\xi<\kappa\rangle$, and for every\/ $\xi<\kappa$ we have
that\/ {\rm $\bfone\forces_{P_\xi}``
\dot Q_\xi$ is hemi-proper.''}
Then\/ $P_\kappa$ is hemi-proper and hence does not collapse\/ $\omega_1$.

Proof. By lemmas 14 and 15.

\section{Namba forcing and the theorem of Ben-David}

In this section we give an application of theorem 16.  The application
is to a theorem of Ben-David which is sketched in [12, theorem XI.1.7].
Because new reals are added at certain
successor stages of the iteration, and possibly even at limit stages
also, the preservation theorems of [12, chapter XI] are not sufficient
to establish this result. The relevant preservation theorem
is [12, XII.3.6, page 408].

Namba forcing, like Prikry forcing, adds no reals
and changes the cofinality of some regular cardinal
to $\omega$. Whereas Prikry forcing collapses a measurable,
Namba forcing collapses $\omega_2$. Hence iterated Namba forcing
uses less extravagant large cardinal assumptions than
iterated Prikry forcing, and indeed for some
assertions proved consistent by Shelah using iterated Prikry forcing,
Shelah later obtained equiconsistency results by using
Namba forcing instead (see [12, chapters X and XI]).

Let us recall the definition of Namba forcing [9].
Let $S=\bigcup\{{}^n\omega_2\,\colon n\in\omega\}$ be the
tree of all finite sequences of ordinals of cardinality at most
$\aleph_1$.

\proclaim Definition 17. Namba forcing is the poset\/
$\{T\subseteq S\,\colon T $ {\rm is non-empty
and $(\forall\eta\in T)(\forall\tau\subseteq\eta)(\tau\in T)$
and $(\forall\eta\in T)($there are $\aleph_2$-many
$\tau\in T$ such that $\tau\supseteq\eta)\}$,}
ordered by inclusion.

Elements of Namba forcing are called {\it perfect subtrees} of $S$.

\proclaim Definition 18. Suppose\/ $P$ is Namba forcing and\/
$G$ is a\/ $V$-generic filter over\/ $P$. Then the\/
{\rm generic object for} $P$ is\/ $\{\alpha\in(\omega_2)^V\,\colon
(\forall T\in G)(\exists\eta\in T)(\alpha\in {\rm range}(\eta))\}$. 

Clearly the generic object is a countable set cofinal
in $(\omega_2)^V$.

\proclaim Theorem 19. 
Assume CH. Then Namba forcing is hemi-proper.

Proof:
The proof follows Namba's proof that the forcing
adds no reals (see, e.g., 
[3, pp.~289--291]).

Suppose $M\prec H_\lambda$ and $N\prec H_{\lambda}$ and $ q$ 
are as in definition 10. 
Take $f\in M$ such that $f$ is a one-to-one map from
$\omega_1$ onto ${}^\omega 2$, and take $g\in M$ a one-to-one
map from $\omega$ onto $\omega\times\omega$.

Work in $N$.
Let $\langle \sigma_n\,\colon n\in\omega\rangle$ list $M^P$. Take
$\langle\xi_n\,\colon n\in\omega\rangle$ a sequence from $M^P$ such that 
$\bfone\forces_P``$if $\sigma_n\in \omega_1$ then
$\xi_{g(n,i)}=f(\sigma_n)(i)$, and in any case $\xi_{g(n,i)}\in 2$''
for every $n$ and $i$ in $\omega$.
For each $n\in\omega$ construct
${\cal Y}_n$ and ${\overline t}_n$ such that
${\cal Y}_n=\langle T_s\,\colon s\in{}^n\omega_2
\rangle$ is a sequence 
of elements of $P$ and
${\overline t}_n=\langle t_s\,\colon s\in {}^n
\omega_2\rangle$ is a sequence of pairwise incompatible elements of
${}^{<\omega}\omega_2$
such that $(\forall s)($every element of $T_s$ is
comparable with $t_s)$ and $(\exists\alpha_s\in 2)(T_s\forces
``\xi_n=\check\alpha_s$'') and $(\forall s'\subseteq s)
(T_s\subseteq T_{s'})$. Because of the final clause
requiring $T_s$ to be stronger than
$T_{s'}$ whenever $s'$ is an initial segment of $s$,
the construction actually proceeds by recursion on $n\in\omega$.
In the base case, of course, we require that
$T_{<>}\leq q$.

For each $n\in\omega$ define ${\cal T}_n$ such that
${\cal T}_n$ is a function with domain
${}^n 2$ such that $(\forall\overline\beta\in
{}^n 2)({\cal T}_n(\overline\beta)=
\bigcup\{T_s\,\colon s\in {}^n\omega_2$ and
$\overline\beta(i)=\alpha_{s\restr i}$ for all $i\leq n\}$.
Take ${\cal T}'$ such that
${\cal T}'$ is a function with domain
${}^\omega 2$ such that $(\forall g\in {}^\omega 2)
({\cal T}'(g) =\bigcap\{{\cal T}_n(g\restr n)\,\colon
n\in\omega\}$).

Claim 1. $(\exists g\in{}^\omega 2)({\cal T}'(g)$
contains a perfect subtree).

Proof: Suppose not. 
For each
$g\in{}^\omega 2$ let $T_0(g)={\cal T}'(g)$ and
for each $\alpha$ let $T_{\alpha+1}(g)=\{t\in T_\alpha(g)\,\colon
t$ has $\aleph_2$-many extensions in $T_\alpha(g)\}$, and
for each limit $\alpha$ let $T_\alpha(g)=\bigcap\{T_\beta(g)\,\colon
\beta<\alpha\}$.
For each $t\in{\cal T}'(g) $ let $h_g(t)$ be the least $\alpha$
such that $t\notin T_\alpha(g)$. This is defined for
every $t\in{\cal T}'(g)$ because otherwise
$\{t\in{\cal T}'(g)\,\colon h_g(t)$ is
undefined$\}$ would be a perfect subtree of ${\cal T}'(g)$,
contrary to assumption.  The relevant facts about the
functions $h_g$ are that $h_g(s)\geq h_g(t)$ whenever $
s\subseteq t\in{\cal T}'(g)$ and for each $t\in{\cal T}'(g)$ there
are at most $\aleph_1$-many extensions $s\supseteq t$ in
${\cal T}'(g)$ such that $h_g(s)=h_g(t)$.  By recursion,
build $s_0\subseteq s_1\subseteq\cdots\subseteq s_n\subseteq\cdots
$ such that for every $n$ we have $s_n\in{}^n\omega_2$ and
$h_g(t_{s_{n+1}})<h_g(t_{s_n})$ for every $g$
for which $h_g(t_{s_{n+1}})$ is defined.
This is possible because, by CH, there are only
$\aleph_1$-many $g$'s in all. Let $g_0\in{}^\omega 2$
be defined by $g_0(n)=\alpha_{s_n}$ for all $n$.
Then $h_{g_0}(t_{s_n})$ is defined for all $n$
because $T_{s_n}\subseteq{\cal T}_{n+1}(\langle\alpha_{s_i}\,\colon
i\leq n\rangle)\subseteq{\cal T}'(g_0)$.
Thus we have a decreasing sequence of ordinals
$h_{g_0}(t_{s_0})>h_{g_0}(t_{s_1})>\cdots>h_{g_0}(t_{s_n})
>\cdots$, a contradiction which establishes the claim.

Still working in $N$,
fix $g$ to witness the claim and take $q'\in P$ such that
${\cal T}'(g)\supseteq q'$.
Fix $n\in\omega$.
Let $\tau_n=\alpha_{g\restr n}$. We have
$q'\subseteq{\cal T}'( g)\subseteq
{\cal T}_n( g\restr n)\leq q$ and $X=_{\rm def}
\{T_s\,\colon s\in{}^n\omega_2\}$ is pre-dense
below ${\cal T}_n( g\restr n)$ and $(\forall s\in X)
( s\forces``\check\tau_n= g(n-1)=\xi_n$'').

We have $\langle\tau_n\,\colon n\in\omega
\rangle\in N$. Hence for each $n\in\omega$ we have
$\tau^*_n=_{\rm def}\langle\tau_{g(n,i)}\,\colon i\in\omega\rangle\in
N$. Thus $q'\forces``$if $\sigma_n\in\omega_1$ then
$\sigma_n=f^{-1}(\tau^*_n)\in N$.''
We have shown $(\forall\sigma\in M^P)(q'\forces``$if $\sigma\in
\omega_1$ then $\sigma\in N$''). Suppose $q'\notforces``\omega_1
\cap M[G_P]\subseteq N$.'' Then we may take $x\in V^P$ to
be a witness. In particular $q'\forces``x\in
M[G_P]$,'' so
we may take $q^*\leq q'$ and $y\in M^P$ such that
$q^*\forces``x=y$.'' Then $q^*\forces``y\in\omega_1$ and
$y\notin N$,'' contrary to what we have already established.
This contradiction shows that $q'\forces``\omega_1\cap M[G_P]
\subseteq N$.'' Hence $q\notforces``\omega_1\cap M[G_P]>\omega_1\cap N$,''
and we are done.

The following is due to Ben-David [12, theorem XI.1.7].

\proclaim Theorem 20. Suppose ZFC$\,+$``there is an inaccessible
cardinal'' is consistent.  Then so is
ZFC$\,+$``there is no cardinal-preserving extension of
the universe in which there is $A\subseteq\omega_2$ such that
$L[A]\models CH$ and $\omega_2^{L[A]}=\omega_2$.''

Proof:
Let $Col(\kappa,\lambda,\theta)$ be the poset which collapses
$\lambda$ to $\kappa$ using conditions of size less than $\theta$.
That is, $Col(\kappa,\lambda,\theta)=\{f\,\colon{\rm dom}(f)\subseteq
\kappa$ and ${\rm range}(f)\subseteq\lambda$ and $\vert f\vert<\theta\}$,
ordered by reverse inclusion.

Begin with a ground model which satisfies $V=L$.
For $\eta$ a limit ordinal or zero, let $\dot Q_\eta$ be Cohen forcing,
and let $Z_\eta$ be (a name for) the corresponding generic subset
of $\omega$. For positive integer $j$ let $\dot Q_{\eta+j}$
be Namba forcing iff $j\in Z_\eta$ and let $\dot Q_{\eta+j}$
be $Col(\omega_1,\omega_2,\omega_1)$ iff $j\notin Z_\eta$.
Then $L[G_{P_\kappa}]\models``(\exists\aleph_2$-many distinct
reals $r\subseteq\omega$ such that there is $\lambda$ which is
a cardinal in $L$ and such that
$(\forall j\in\omega)(r(j)=0$ iff $\cf^{L[G_{P_\kappa}]}(\lambda^{+j})^L
=\omega$)).'' To clarify, $(\lambda^{+j})^L$ is the $j^{th}$ successor
of $\lambda$ as computed in $L$, and its cofinality is to be
computed in $L[G_{P_\kappa}]$. Now suppose, towards
a contradiction, that there is a cardinal-preserving extension
$V\supseteq L[G_{P_\kappa}]$ and $A\in V$ such that
$A\subseteq\omega_2^V=\kappa$ and $L[A]\models CH$ and
$\omega_2^{L[A]}=\omega_2^V=\kappa$. 
Let $X=L_\theta[\{\delta<\omega_2^V\,\colon\cf^{L[G_{P_\kappa}]}
(\delta)=\omega$ and $\delta$ is a cardinal of $L\}]$ where
$\theta$ is some big cardinal (say, $\kappa^+$).
Because $L$ and $L[G_{P_\kappa}]$ are
both subuniverses of $V$
we know that $X\in V$.
But $X = L_\theta[\{\delta<\omega_2^V\,\colon\cf^V(\delta)=\omega$
and $\delta$ is a cardinal of $L\}]$ because $V$ is a
cardinal preserving extension. Because $\omega_2^{L[A]}=\omega_2^V$
we have that $\cf^{L[A]}(\delta)=\cf^V(\delta)$ for all $\delta<\omega_2^V$.
Hence in $L[A]$ we have
$\{\delta<\omega_2^{L[A]}\,\colon\cf^{L[A]}(\delta)=\omega\}$
and we have $\{\delta<\omega_2^{L[A]}\,\colon\delta$ is a
cardinal of $L\}$. From these two sets we may recover
$\aleph_2^V$-many reals (namely, $\{Z_\eta\,\colon\eta$ is a limit
ordinal less than $\kappa\}$). Hence $L[A]\not\models CH$.
This contradiction establishes the theorem.

\section{$P[S]$ is hemi-proper}

We show that $P[S]$ is hemi-proper. Combined with the argument of
[12, section XI.7], this gives a new proof of Shelah's answer
to the problem of H. Friedman [2]. 

\proclaim Definition 21. Suppose\/ $S$ is a stationary
subset of\/ $S^2_0=_{\rm def}\{\alpha<\omega_2\,\colon\cf(\alpha)=\omega\}$.
Then\/ $P[S]$ is the poset\/ {\rm $\{f\,\colon(\exists\beta<\omega_1)
(\dom(f)=\beta+1$ and ${\rm range}(f)\subseteq S$ and
$f$ is continuous increasing$)\}$,} ordered by reverse inclusion
(thus the ranges of the conditions are ordered by reverse
end-extension).

\proclaim Theorem 22. Suppose $S\subseteq S^2_0$ is stationary.
Then $P[S]$ is hemi-proper.

Proof.  Suppose $\lambda$, $M$, $N$, and $q$ are as in definition 10.
Take $\mu<\lambda$ regular such that
the power set of $P[S]$ is in $H_\mu$ and $\mu\in M$ and $\mu^+<\lambda$.
Take $\langle M_i\,\colon i\in\omega_2\rangle\in N$ a continuous
tower of elementary subtructures of $H_{\mu^+}$ each of cardinality $\aleph_1$,
such that $M\cap H_\mu\in M_0$. Take $\delta\in S\cap N$ such that
$\sup(\omega_2\cap M_\delta)=\delta$. Take $\langle\delta_n\,\colon
n\in\omega\rangle$ an increasing
sequence from $N$
cofinal in $\delta$. Let $\langle D_n\,\colon n\in\omega\rangle$
enumerate all dense open sets of $P[S]$ which are in $N\cap M_\delta$.
We may suppose 
that $(\forall n\in\omega)(D_n\in M_{\delta_n})$.
Build $q\geq q_0\geq q_1\geq\cdots$ such that $q_n\in
D_n\cap M_{\delta_n}\cap 
N$ and $\sup(q_{n+1})\geq\sup(\omega_2\cap M_{\delta_n})$.
Then let $r=\bigcup\{q_n\,\colon n\in\omega\}\cup
\{\delta\}$. 
We have $r\forces``OR\cap M[G_{P[S]}]\subseteq
N$'' where $OR$ is the class of ordinals. We are done.

\proclaim Theorem 23 {\rm (Shelah)}. Suppose a Mahlo cardinal is consistent.
Then so is the statement: Whenever\/ $S$ is a stationary subset
of\/ $S^2_0$ then\/ $S$ contains an uncountable sequentially closed subset.

Proof: Combine the argument of [12, section XI.7] with our
theorems 16 and 22.

\section{Pseudo-completeness}

In this section, we introduce the notion of
$\mu$-pseudo-completeness for $\mu$ a cardinal
$(\mu = 2$ or $\mu = {\rm OR}$ is allowed). This is a generalization
of $\mu$-completeness, and similar to [12, section X.3].
We then define $\mu$-pseudo-completeness for $\dot Q$
relative to a poset $P$, in a manner reminiscent of
the ``not adding reals'' theorem of [10].
We apply these results
to iterated Namba forcing in the following section and
to $P[S]$ in a later section.

\proclaim Definition 24. Suppose $\mu$ is a regular cardinal
(we allow $\mu=2$ or $\mu=OR$) and\/ $P$ is a poset and\/ $\dot Q$
is a $P$-name for a poset.  We say that\/ $\dot Q$
is\/ {\rm $\mu$-pseudo-complete relative to $P$} iff whenever
$\lambda$ is a sufficiently large regular cardinal and\/
$M$ is a countable elementary substructure of\/
$H_\lambda$ and\/ $N$ is a countable
elementary substructure of\/ $H_{\lambda^+}$
and $\{P*\dot Q,\mu\}\subseteq M\in N$
and\/ $p\in P$ and $\dot q\in N^P$ and
$p\forces``\dot q\in\dot Q\cap M[G_P]$,'' then there is
$\dot r\in V^P$ such that $p\forces``\dot r\leq\dot q$''
and whenever $\sigma\in N^{P*\dot Q}$ and\/ {\rm $\bfone
\forces``\sigma\in \mu\cap
M[G_{P*\dot Q}]$''} then there is
$\tau\in N^P$ such that\/    {\rm
$\bfone\forces_P``\tau\in M[G_P]$''} and\/ {\rm $p\forces``
\dot r\forces`
\sigma=\check\tau$.'\thinspace''}

Remark: If $\mu={\rm OR}$ we waive $\mu\in M$.
We shall omit stating this exception in the sequel.

\proclaim Definition 25. Suppose\/ $\langle P_\xi\,\colon
\xi\leq\kappa\rangle$ is a countable support iteration.
We say that\/ $P_\kappa$ is\/ {\rm strictly $\mu$-pseudo-complete}
iff whenever\/ $\eta<\kappa$ and\/ $\lambda$ is a sufficiently
large regular cardinal and\/ $M$ is a
countable elementary
substructure of\/ $H_{\lambda}$ and\/
$N$ is a countable elementary substructure of\/ $H_{\lambda^+}$
and\/ $\{P_\kappa,\mu,\eta\}\subseteq
M\in N$  and\/  $p\in P_\eta$ and\/
 $\dot q\in
N^{P_\eta}$ and\/  {\rm
$p\forces_{P_\eta}``\dot q\in\dot P^{M}_{\eta,\kappa}$,''}
then there is\/ $r\in P_\kappa$ such that\/ $r\restr\eta= p$ and\/
$p\forces``r\restr[\eta,\kappa)\leq\dot q$'' and\/
$\supt(r)\subseteq\eta\cup M$
and whenever $\sigma\in N^{P_\kappa}$ and $\bfone\forces``\sigma\in
\mu\cap M[G_{P_\kappa}]$'' then there is\/
$\tau\in N^{P_\eta}$ such that\/  {\rm
$\bfone\forces_P``\tau\in\mu\cap M[G_{P_\eta}]$''} and\/ {\rm $p\forces``
r\restr[\eta,\kappa)\forces`
\sigma=\check\tau$.'\thinspace''}

The reason the following lemma does not contradict [5]
 is that definition 25 must hold even when $p\notin N$.

\proclaim Lemma 26. Suppose\/ $\langle P_\xi\,\colon\xi\leq\kappa
\rangle$ is a countable support iteration based on\/
$\langle \dot Q_\xi\,\colon\xi<\kappa\rangle$, and
suppose\/ $P_\xi$ is strictly\/ $\mu$-pseudo-complete
for every\/ $\xi<\kappa$, and suppose that if\/ $\kappa=\gamma+1$
then\/ $\dot Q_\gamma$ is\/ $\mu$-pseudo-complete relative to\/
$P_\gamma$. Then\/ $P_\kappa$ is strictly\/ $\mu$-pseudo-complete.

Proof: Suppose $\lambda$, $\eta$, $M$, $N$, $\dot q$, and
$p$ are as in definition 25. 

Case 1: $\gamma+1=\kappa$ for some $\gamma$.

We may
take $r_0\in P_\gamma$ such that
$r_0\restr\eta=p$ and $p\forces``r_0\restr[\eta,\gamma)\leq
\dot q\restr\gamma$'' and $\supt(r_0)\subseteq
\eta\cup M$ and whenever
$\sigma\in N^{P_\gamma}$ there is
$\tau\in N^{P_\eta}$ such that if
$\bfone\forces``\sigma\in\mu\cap
M[G_{P_\gamma}]$'' then $\bfone\forces ``\tau\in \mu\cap M[G_{P_\eta}]$''
and $p\forces``
r_0\restr[\eta,\gamma)\forces`\sigma=\check\tau$.'\thinspace''
Take $\dot s\in V^{P_\gamma}$
such that $r_0\forces_{P_\gamma}``\dot s\leq\dot q(\gamma)$''
and for every $\sigma\in N^{P_\kappa}$ such that if
$\bfone\forces``\sigma\in\mu\cap M[G_{P_\kappa}]
$'' there is
$\tau\in N^{P_\gamma}$ such that $\bfone\forces``\tau\in \mu\cap
M[G_{P_\gamma}]$'' and $r_0\forces``\dot s\forces`
\sigma = \check\tau$.'\thinspace''
Take $r\in P_\kappa$ such that $r\restr\gamma=r_0$ and
$r(\gamma)=\dot s$. Then $r$ is as required.

Case 2: $\kappa$ is a limit ordinal.

Let $\langle\sigma_n\,\colon n\in\omega\rangle$ list
$\{\sigma\in N^{P_\kappa}\,\colon\bfone\forces``\sigma\in\mu\cap
M[G_{P_\kappa}]$''$\}$.
 Let $\alpha=\sup(\kappa\cap M)$, and let
 $\langle \alpha_n\,\colon n\in\omega\rangle$ be an increasing
sequence of ordinals from $\alpha\cap M$ cofinal in $\alpha$
such that $\alpha_0=\eta$.

Build $\langle p_n,\dot q_n,\tau_n\,\colon
n\in\omega\rangle$
such that $p_0=p$ and
$\dot q_0=\dot q$ and
both of the following:

($i$) $p_{n}\forces_{P_{\alpha_n}}
``\dot q_{n+1}\leq\dot q_n\restr[\alpha_{n},
\kappa)$ and $\dot q_{n+1}\in\dot P_{\alpha_{n},\kappa}^{M}$
and $\dot q_{n+1}\forces`
\sigma_n=\check\tau_n$'\thinspace''
and $\tau_n\in N^{P_{\alpha_{n}}}$
and $\dot q_{n+1}\in N^{P_{\alpha_n}}$

($ii$)  $p_n\in P_{\alpha_n}$ and $p_{n+1}\restr\alpha_n=p_n$
 and $\supt(p_{n+1})\subseteq\eta\cup M$  and
 $p_n\forces``p_{n+1}\restr[\alpha_n,\alpha_{n+1})\leq
 \dot q_{n+1}\restr\alpha_{n+1}$''
and whenever $\sigma\in N^{P_{\alpha_{n+1}}}$ and
$\bfone\forces``\sigma\in \mu\cap M[G_{P_{\alpha_{n+1}}}]$''
there is $\tau\in
N^{P_{\alpha_n}}$ such that $\bfone\forces``\tau\in
\mu\cap M[G_{P_{\alpha_n}}]$'' and $p_{n+1}
\forces``
\sigma=\tau$''

We may choose $\dot q_{n+1}$ and $\tau_n$ as in $(i)$ by lemma 2.
We may choose $p_{n+1}$ as in $(ii)$ by the fact that
$P_{\alpha_{n+1}}$ is strictly $\mu$-pseudo-complete.

Now take $r\in P_\kappa$
such that $\supt(r)\subseteq\alpha$ and
for all $n\in\omega$ we have $r\restr\alpha_n=p_n$.
The lemma is established.

\proclaim Theorem 27. Suppose\/ $\langle P_\xi\,\colon\xi\leq\kappa\rangle$
is a countable support iteration based on\/ $\langle\dot Q_\xi\,\colon
\xi<\kappa\rangle$, and for every\/ $\xi<\kappa$ we have
that\/ $\dot Q_\xi$ is
 $\mu$-pseudo-complete relative to\/ $P_\xi$.
Then\/ $P_\kappa$ does not add any elements of\/ ${}^\omega\mu$.

Proof. By lemma 26 we have that
$P_\kappa$ is strictly $\mu$-pseudo-complete.
Take $\eta=0$ in definition 25.

\section{Namba forcing is $\mu$-pseudo-complete ($\mu\in\{2,\omega\})$}

We now investigate applications of theorem 27.
It is easy to show that if $\mu<\mu^*$ are cardinals and
$\bfone\forces_P``\mu^*$ is measurable and $\dot Q$ 
is Prikry forcing on $\mu^*$ (more exactly, on
some fixed normal measure over $\mu^*)$'' then $\dot Q$ is
$\mu$-pseudo-complete relative to $P$.

\proclaim Lemma 28. Suppose $P$ is a poset which does not add reals
and $\bfone\forces_P``\dot Q$ {\rm is Namba forcing.''}
Suppose also CH holds in the ground model and
$\mu\in\{2,\omega\}$. Then $\dot Q$ is
$\mu$-pseudo-complete relative to $P$.

Proof: As in theorem 19, we modify
slightly the argument of [9], [3, pp.~289--291].
Suppose $M\prec H_\lambda$ and $N\prec H_{\lambda^+}$
and $p\in P$ and $\dot q\in N^P$
are as in definition 24. Let
$\langle\sigma_n\,\colon n\in\omega\rangle$ list $\{\sigma\in
N^{P*\dot Q}\,\colon\bfone\forces``\sigma\in\mu$''$\}$.
Take $\dot\beta\in M^P$ such that $\bfone\forces_P``\dot \beta
=(\omega_2)^{V[G_P]}$'' (notice that we allow the possibility
of $p_1\in P$ and $p_2\in P$ and $\beta_1\ne\beta_2$ such that
$p_1\forces``\beta_1=(\omega_2)^{V[G_P]}$'' and
$p_2\forces``\beta_2=(\omega_2)^{V[G_P]}$'').
By recursion on $n\in\omega$ construct
${\cal Y}_n\in N^P$ and ${\overline t}_n\in N^P$ such that
$\bfone\forces_P``$if $\dot q\in
\dot Q\cap M[G_P]$ then ${\cal Y}_n=\langle T_s\,\colon s\in{}^n\dot\beta
\rangle$ is a sequence 
of elements of $\dot Q$ and
${\overline t}_n=\langle t_s\,\colon s\in {}^n
\dot\beta\rangle$ is a sequence of pairwise incompatible elements of
${}^{<\omega}\dot\beta$
such that $(\forall s)($every element of $T_s$ is
comparable with $t_s)$ and $(\exists\alpha_s\in\mu)(T_s\forces
`\sigma_n\in\mu$ implies $
\sigma_n=\check\alpha_s$') and $(\forall s'\subseteq s)
(T_s\subseteq T_{s'})$;'' also
$\bfone\forces_P``$if $\dot q\in\dot Q\cap M[G_P]$
then $T_{<>}\leq\dot q$.''
Furthermore, we may assume that $\overline\alpha_n$ is a name
such that $\bfone\forces_P``\overline\alpha_n=
\langle\alpha_s\,\colon s\in {}^n\dot\beta\rangle$ as in the
preceding construction'' and $\overline\alpha_n\in N^P$.

For each $n\in\omega$ define ${\cal T}_n\in N^P$ such that
$\bfone\forces_P``{\cal T}_n$ is a function with domain
${}^n\mu$ such that $(\forall\overline\beta\in
{}^n\mu)({\cal T}_n(\overline\beta)=
\bigcup\{T_s\,\colon s\in {}^n\dot\beta$ and
$\overline\beta(i)=\alpha_{s\restr i}$ for all $i\leq n\}$.''
We have stayed within the confines of $N^P$ as long as possible; now
we start making names which are in $V^P$.
Take ${\cal T}'\in V^P$ such that
$\bfone\forces_P``{\cal T}'$ is a function with domain
${}^\omega\mu$ such that $(\forall g\in {}^\omega\mu)
({\cal T}'(g) =\bigcap\{{\cal T}_n(g\restr n)\,\colon
n\in\omega\}$).''

Claim: $(\exists g\in{}^\omega\mu)(\bfone\forces``$if
$\dot q\in\dot Q\cap M[G_P]$ then
${\cal T}'(\check g)$
contains a perfect subtree'').

Proof: Suppose not. We have $(\forall g\in{}^\omega\mu)
(\bfone\notforces``$if $\dot q\in\dot Q\cap
M[G_P]$ then ${\cal T}'(\check g)$ contains a perfect
subtree''). Because $P$ adds no reals, we may take
$q_2\in P$ such that $q_2\forces``\dot q\in\dot Q\cap M[G_P]$
and $(\forall g\in
{}^\omega\mu)({\cal T}'(g)$ does not contain a perfect
subtree).''
Take $G_P$ to be $V$-generic over $P$ such that $q_2\in G_P$.
Using theorem 19 (claim 1) in the model $V[G_P]$ we obtain
a contradiction (literally so if $\mu=2$, but if
$\mu=\omega$ then rewrite the proof of that claim with $\omega$
replacing 2).

Fix $g$ to witness the claim and take $\dot q_1\in V^P$ such that
$\bfone\forces_P``$if $\dot q\in
\dot Q\cap M[G_P]$ then ${\cal T}'(\check g)\supseteq \dot q_1$
and $\dot q_1\in \dot Q$.''
Given $\sigma\in N^{P*\dot Q}$ such that
$\bfone\forces``\sigma\in\mu$,''
take $n$ such that $\sigma=\sigma_{n-1}$.
We seek $\tau\in N^P$ such that $p\forces``\dot q_1\forces`
\sigma=\check\tau$.'\thinspace''
Take $\tau\in N^P$ such that $\bfone\forces_P``\tau=\alpha_{g\restr n}$.''
Although $g$ need not be in $N$, certainly $g\restr n\in N$,
so there is no problem in choosing such a $\tau$.  We have
$p\forces``\dot q_1\subseteq{\cal T}'(\check g)\subseteq
{\cal T}_n(\check g\restr n)$ and $X=_{\rm def}
\{T_s\,\colon s\in{}^n\dot\beta\}$ is pre-dense
below ${\cal T}_n(\check g\restr n)$ and $(\forall\dot s\in X)
(\dot s\forces`\check\tau=\check g(n-1)=\sigma$').''
The lemma is established.

\section{Applications}

In this section we give several applications of
iterated Namba forcing.  All
of these are taken from [12, chapter XI] and are included merely for
the sake of completeness.
Our first application is [12, theorem XI.1.5].  Our proof is
the same as the one given by Shelah, except of course that he 
uses revised countable support iterations and he
intersperses cardinal collapses among the Namba forcings, as
required by his
preservation theorem.
The converse of the theorem had earlier been proved by Avraham,
so that this is an equiconsistency result.

\proclaim Theorem 29 {\rm (Shelah)}.
Suppose ZFC$+$``there exists an inaccessible
cardinal'' is consistent.  Then 
so is ZFC$+$GCH$+(\forall X\subseteq\omega_1)
(\exists Y\in[\omega_2]^\omega)(Y\notin L[X])$.

Proof: Start with a ground model of ZFC$+$GCH$+\kappa$ is inaccessible.
Form a countable support
iteration of length $\kappa$ such that each
constituent poset is Namba forcing.
Clearly $\omega_1^V=\omega_1^{V[G_{P_\kappa}]}$ and
$\kappa=\omega_2^{V[G_{P_\kappa}]}$, and GCH holds in $V^{P_\kappa}$.
Suppose $X\in V[G_{P_\kappa}]$ and $V[G_{P_\kappa}]\models``\vert
X\vert=\aleph_1$.''
Take $\alpha<\kappa$ such that $X\in V[G_{P_\alpha}]$.
In $V[G_{P_\kappa}]$ let $Y$ be the generic object for
$\dot Q_{\alpha}$. Then we have $Y\notin V[G_{P_\alpha}]$,
but $L[X]\subseteq V[G_{P_\alpha}]$. Thus $Y$ exemplifies
what is required, and the theorem is established.

The following theorem is [12, theorem XI.1.6]. Again, we add nothing
new (we include it for expository purposes) except that
our iteration of Namba forcings uses countable
support iteration and does not need the other $\sigma$-closed
cardinal collapses.  Avraham has proved the converse, so this
theorem is an equiconsistency result.

\proclaim Theorem 30 {\rm (Shelah)}.
Suppose ZFC$+$``there is a Mahlo cardinal''
is consistent. Then so is ZFC$+$GCH$+(\forall A\subseteq
\omega_1)(\exists\delta>\omega)(\cf(\delta)=\omega$
{\rm and $\delta$ is a regular cardinal in $L[A\cap \delta])$.}

Proof: Start with a ground model satisfying ZFC$+$GCH$+\kappa$
is Mahlo. Form a countable support
iteration of length $\kappa$ using
Namba forcings.  Given $A\in V^{P_\kappa}$ such that
$\bfone\forces``A\subseteq\omega_1$,''
take $C\subseteq\kappa$ closed unbounded such that
whenever $i<j$ are both in $C$ then
$A\cap i\in V^{P_j}$.   Take $C'$ a closed unbounded set
such that every element of $C'$ is a limit point of $C$, and
such that $A\in V^{P_\eta}$ where $\eta={\rm min}(C')$.
Then we have that $(\forall\delta\in C')(A\cap\delta\in
V^{P_\delta})$. Since $\kappa$ is Mahlo, we may
take $\lambda\in C'$ such that $\lambda$ is inaccessible.
We have that $P_\lambda$ has $\lambda$-c.c., hence
$\omega_2^{P_\lambda}=\lambda$. Let $Y$ be the generic
object for $\dot Q_\lambda$. Then
$V^{P_\kappa}\models``Y$ is cofinal in $\lambda$ hence
$\cf(\lambda)=\omega$.'' Yet $\lambda $ is regular
in $L[A\cap\lambda]$ because $L[A\cap\lambda]\subseteq V^{P_\lambda}$
and $V^{P_\lambda}\models``\lambda$ is regular.''
The theorem is established.

We now give another application of iterated Namba forcing
from Shelah's book. Once again, this is  for expository purposes;
however, we use an argument from
[4] rather than the argument used in [12].

\proclaim Definition 31. Suppose ${\cal F}$ is a filter over $\kappa$.
Then the poset $PP({\cal F})$ is $\{X\subseteq\kappa\,\colon
\kappa-X\notin{\cal F}\}$, ordered by inclusion.

\proclaim Definition 32.  Suppose ${\cal F}$ is a filter over $\kappa$.
We say ${\cal F}$ is {\rm precipitous} iff
$\bfone\forces_{PP({\cal F})}``V^\kappa/E$ {\rm is
well-founded where $E$ is $PP({\cal F})$-generic over $V$.''}

\proclaim Lemma 33. Suppose\/
 $\langle P_\alpha\,\colon\alpha\leq\kappa\rangle$
is a forcing iteration such that $(\forall\alpha<\kappa)(\vert P_\alpha\vert
<\kappa)$ and $(\forall p\in P_\kappa)({\rm supt}(p)$
{\rm is bounded below $\kappa$)} and $U$ is a normal
measure on $\kappa$ and $\bfone\forces_{P_\kappa}``{\cal F}$
{\rm is the filter generated by $\check U$.''}
Then $\bfone\forces_{P_\kappa}``{\cal F}$ {\rm is precipitous.''}

Remark: [3, page 592] attributes the following argument,
which is the first half of [4, theorem 3], to Mitchell,
but there Levy forcing is used in place of the $\kappa$-length
iteration, and the filter is constructed over $\omega_1$ instead
of $\omega_2$.
A different argument, based on [8], is used by
Shelah [12, theorem XI.1.7].

Proof: Let $j\,\colon V\mapsto V^\kappa/U$ be the canonical embedding
and let $G=G_{P_\kappa}$ be the canonical name for the
generic filter on $P_\kappa$. Because $(\forall\alpha<\kappa)
(\vert P_\alpha\vert<\kappa)$ we have that $j(P_\kappa)=
P_{j(\kappa)}=P_\kappa*\dot P_{\kappa,j(\kappa)}$. Let
$G^*=G_{j(P_\kappa)}$ be the canonical name for the
generic filter on $j(P_\kappa)$. In $V^{j(P_\kappa)}$ we have
$G=G^*\cap P_\kappa$. Let $j^*\in V^{P_\kappa}$
be a name for the elementary embedding from $V[G]$ into
$(V^\kappa/U)[G^*]$ which extends $j$ and such that
$j^*(G)=G^*$. Actually of course $j^*$ is a propewr class
and therefore not literally an element of $V^{P_\kappa}$
but his need not concern us unduly. We shall use
standard facts about $j^*$ which can be found in [7],
[3, chapters 36 and 37].

Claim 1.
Suppose $y\in V^{P_\kappa}$ and  $p\in P_\kappa$.
Then $(\exists X\in U)(p\forces_{P_\kappa}``y\supseteq\check X$'')
iff $j(p)\forces_{j(P\kappa)}``\kappa\in j^*(y)$.''

Proof: Suppose first that $j(p)\forces``\kappa\in j^*(y)$.''
Let $X=\{\alpha<\kappa\,\colon p\forces``\alpha\in y$''$\}$.
By normality of $U$ our hypothesis implies $X\in U$ (see [7]).
But $p\forces``\check X\subseteq y$'' by the definition of $X$.
In the converse direction, suppose $Z\in U$ and $p\forces``\check Z
\subseteq y$'' and $q\leq j(p)$.
It suffices to find $q^*\leq q$ such that $q^*\forces``\kappa\in
j^*(y)$.'' Let $q'=q\restr\kappa\in P_\kappa$.
We have $\bfone\forces_{j(P_\kappa)}``\kappa\in
j^*(\check Z)$'' by normality of $U$. Also,
$j(q')\forces_{j(P_\kappa)}``\kappa\in
j^*(y)$,'' because
$q'\forces``\check Z\subseteq y$.'' Take $q^*$ such that
$q^*\restr\kappa= q'$ and $q^*\restr[\kappa,j(\kappa))=
q\restr[\kappa,j(\kappa))$. Because
$\supt(q')\subseteq\gamma$ for some $\gamma<\kappa$
we have $\supt(j(q'))\subseteq\gamma$ and hence
$q^*\in j(P_\kappa)$  and $q^*\leq q$ and $q^*\leq j(q')$.
We have $q^*\forces``\kappa\in j^*(\check Z)\subseteq
j^*(y)$.'' The claim is established.

Suppose $p\in P_\kappa$ and $p\forces``x\in PP({\cal F})$.''
We show there is $q\leq p$ and $D\in V^{P_\kappa}$ such that
$q\forces``x\in D$ and $D$ is $PP({\cal F})$-generic and
$V[G]^\kappa/D$ is wellfounded,'' which suffices to
establish the lemma. (Remark: It suffices to find such a
$D$ in $V[G][G_{PP({\cal F})}]$. It may therefore
seem surprising that we obtain such a $D$
in $V[G]$. This occurs essentially because in the
ground model $V$ we have $U\in V$ is $PP(U)$-generic over $V$,
and since $P_\kappa$ has $\kappa$-c.c., we have moved to a situation
which is not too far removed from the situation in the ground model.)
By the definition of $PP({\cal F})$ we have
$p\forces_{P_\kappa}``\kappa-x\notin{\cal F}$.''
Hence $p\notforces``\bfone\forces_{\dot P_{\kappa,j(\kappa)}}`
\kappa\in j^*(\kappa-x)$.'\thinspace''
Take $p_1\leq p$ such that $p_1\forces``\bfone\notforces`\kappa\in
j^*(\kappa-x)$.'\thinspace'' Take $q_1\in j(P_\kappa)$
with $q_1\restr\kappa\leq p_1$ and $q_1\forces``\kappa\notin
j^*(\kappa-x)$.'' Take $q_2\leq q_1$ such that $q_2\forces``\kappa
\in j^*(x)$.'' Define $D\in V^{P_\kappa}$
to be a name for $\{y\subseteq\kappa\,\colon
q_2\restr[\kappa,j(\kappa))\forces_{\dot P_{\kappa,j(\kappa)}}``\kappa\in
j^*(y)$''$\}$ and let $q=q_2\restr\kappa$.
 We show that $q$ and $D$ satisfies the requirements.

By claim 1 we have that  $V[G]\models``{\cal F}$ is normal and
$D\supseteq {\cal F}$ because $(\forall X\subseteq\kappa)
(X\in{\cal F}$ iff $\bfone\forces_{\dot P_{\kappa,j(\kappa)}}
`\kappa\in j^*(X)$'),''
using once again the fact that $\supt(j(p'))\subseteq\kappa$
whenever $p'\in P_\kappa$.

Continuing the verification of the requisite properties
of $q$ and $D$, we have by choice of $q_2$ that
$q_2\forces``\kappa\in j^*(x)$.'' Thus $q\forces``x\in D$''
by the definition of $D$.
Furthermore, $q\forces``(V[G])^\kappa/D$ is well-founded
because $\bfone\forces_{\dot P_{\kappa,j(\kappa)}}`(V[G])^\kappa/D$
is isomorphic
to the well-founded structure $(V^\kappa/U)[G^*]$ via the isomorphism
$\pi([f]) = j^*(f)(\kappa)$.'\thinspace''
So it remains to show $q\forces``D$ is $PP({\cal F})$-generic.''

Suppose, towards a contradiction, that $A\in V^{P_\kappa}$ and $q'\leq q$
and $q'\forces``A\subseteq PP({\cal F})$ and $A\cap D=\emptyset$ and
$A$ is dense in $PP({\cal F})$.''
Let $q_3\in P_{j(\kappa)}$ be defined by $q_3\restr\kappa=q'$
and $q_3\restr[\kappa,j(\kappa))=q_2\restr[\kappa,j(\kappa))$.

Claim 2. $q_3\forces``(\forall a\in A)(\kappa\notin j^*(a))$.''

Proof: Suppose instead that $q_4\leq q_3$ and $a\in V^{P_\kappa}$
and $q_4\forces``a\in A$ and $\kappa\in j^*(a)$.''
Then $q_4\restr\kappa\notforces``q_2\restr[\kappa,j(\kappa))\notforces`
\kappa\in j^*(a)$.'\thinspace''
Hence $q_4\restr\kappa\notforces``a\notin  D$,'' a contradiction
to the fact that $q'\forces``A\cap D=\emptyset$''.

Take ${\overline g}=\langle g_\alpha\,\colon\alpha<\kappa\rangle$ such that
$[{\overline g}]_U=q_3 $ and each $g_\alpha$ maps $\kappa$ into
$P_\kappa$. Take $T\in V^{P_\kappa}$ such that
$\bfone\forces``T=\{\alpha<\kappa\,\colon g_\alpha\in G\}$.''
Because $q_3\forces``q_3\in G^*= j^*(G)$'' we have
$q_3\forces``\{\alpha<j(\kappa)\,\colon\alpha\in j(T)\}\in
j(U)$'' and hence $q_3\restr\kappa\notforces``\{\alpha<\kappa\,\colon
\alpha\in T\}\notin U$.'' Hence we may take $q_4\leq q_3$ such that
$q_4\forces``\kappa\in j^*(T)$.'' Hence $q_4\restr\kappa\notforces``
\kappa-T\in{\cal F}$.'' Therefore we may take $q_5\leq
q_4\restr\kappa$ such that $q_5\forces``\kappa-T\notin{\cal F}$.''
Because $q_5\forces``A$ is dense in $PP({\cal F})$,''
we may take $q_6\leq q_5$ and $a\in V^{P_\kappa}$ such that
$q_6\forces``a\in A$ and $\kappa-(a\cap T)\notin{\cal F}$.''
We have $q_6\notforces``\bfone\forces_{\dot P_{\kappa,j(\kappa)}}`
\kappa\in j^*(\kappa-(a\cap T))$.'\thinspace''
Because $a$ and $T$ are in $V^{P_\kappa}$ we have
$q_6\forces``\bfone_{\dot P_{\kappa,j(\kappa)}}$ decides
$`\kappa\in j^*(\kappa-(a\cap T)$'\thinspace''
hence we may take
$q_7\leq q_6$ such that $q_7\forces``\bfone\forces_{\dot P_{\kappa,
j(\kappa}}`\kappa\notin j^*(\kappa-(a\cap T))$.'\thinspace''
Using the fact that $q_4\forces``\kappa\in j^*(T)$'' we have
$q_7\notforces``\kappa\notin j^*(a)$,'' contradicting claim 2.
The lemma is established.

The following is [12, theorem XI.1.7].

\proclaim Theorem 34
{\rm (Shelah)}. Suppose ZFC$+$``there is a measurable cardinal''
is consistent. Then so is ZFC$+$GCH$+$``there is
a normal precipitous filter $D$ over $\omega_2$ such that
$S^2_0=\{\alpha<\omega_2\,\colon\cf(\alpha)=\omega\}\in D$.''

Proof: Begin with a model of ZFC$+$GCH in which $\kappa$ is a measurable
cardinal. Form a countable support
iteration $\langle P_\eta\,\colon\eta\leq
\kappa\rangle$ with $\bfone\forces_{P_\eta}``\dot Q_\eta$ is Namba forcing''
for every $\eta<\kappa$. 
Let $U$ be a normal measure on $\kappa$
and let  $Y=\{\alpha<\kappa\,\colon\alpha$ is strongly inaccessible$\}$.
By lemma 33, $\bfone\forces_{P_\kappa}``{\cal F}$ is
precipitous, where ${\cal F}$ is the filter generated by
$\check U$.'' Thus it suffices to show
$\bfone\forces_{P_\kappa}``\{\alpha<\omega_2^{V[G_{P_\kappa}]}=\kappa
\,\colon\cf(\alpha)=\omega\}\supseteq\check Y$.'' But all that is needed
to see this is to notice that if $\alpha$ is strongly inaccessible
then  we have that $P_\alpha$ has
$\alpha$-c.c., and so $\bfone\forces_{P_\alpha}``\alpha
=\omega_2^{V[G_{P_\alpha}]}$ and so $\dot Q_\alpha$ adds a countable
sequence cofinal in $\alpha$.'' The theorem is established.

\section{$P[S]$ is pseudo-complete}

In this section we establish that $P[S]$ is $\mu$-pseudo-complete
for $\mu=OR$.  For\/ $p\in P[S]$, we set $\hht(p)
={\rm max}(\dom(p))$, i.e., $\dom(p)=\hht(p)+1$.

\proclaim Theorem 35. Suppose\/ $\dot Q$ is a\/ $P$-name and\/ {\rm
$\bfone\forces_P``\dot Q=P[\dot S]$ for some
stationary $\dot S\subseteq S^2_0$ and CH holds.''} Then\/ $\dot Q$ is\/
$\mu$-pseudo-complete relative to
$P$ for every regular\/ $\mu$ including\/
$\mu=OR$.

Proof: Suppose $P*\dot Q\in M\prec H_\lambda$ and
$N\prec H_{\lambda^+}$ with
$M\in N$ both countable, and $(p,\dot q)\in P*\dot Q$ and
$\dot q\in N^P$, as in definition 24.
Let $\langle\sigma_n\,\colon n\in\omega\rangle $
list $\{\sigma\in N^P\,\colon\bfone\forces``\sigma\in M[G_P]$
and $\sigma$ is an ordinal''$\}$.
 Build $\{J_n\,\colon n\in
\omega\}\subseteq N^P$ by recursion on $n$ such that
for each $n\in\omega$ we have $\bfone\forces_P``J_n\subseteq
\dot Q$ and $(\forall x\in J_n)(\exists\eta
<\hht(x))(x\restr\eta\in J_{n-1})$ and
$(\forall t\in J_{n-1})(\vert\{x\in J_n\,\colon
x\restr\dom(t)=t\}\vert=\aleph_2)$ and
$(\forall x\in J_n)(\exists\tau\in\mu)(x\forces_{\dot Q}`
$if $\sigma_n\in\mu$ then $
\sigma_n=
\check\tau$').''
Take $\tau_n\in N^P$ such that
$\bfone\forces_P``(\forall x\in J_n)(x\forces`
$if $\sigma_n\in\mu$ then $
\sigma_n=\check\tau_n$').''
By the fact that $\bfone\forces_P``$CH holds,'' we may
take $C\in V^P$ such that
$\bfone\forces``C\subseteq\omega_2$ is closed
unbounded and whenever $\delta\in C$ and $x\in\bigcup\{J_n\,\colon
n\in\omega\}$ and $\sup({\rm range}(x))<\delta$
then there are $\beta<\delta$ and $y\in\bigcup\{J_n\,\colon
n\in\omega\}$ such that
$y\restr\dom(x)=x$ and $y(\dom(x))=\beta$.''

Take $\delta\in V^P$ such that $\bfone\forces``\delta\in C\cap
\dot S$ and
$\sup({\rm range}(\dot q)<\delta$.'' 
Also suppose $\bfone\forces``\langle\delta_n\,\colon
n\in\omega\rangle$ is an increasing sequence cofinal in $\delta$.''
By recursion on $n\in\omega$ build $\{\dot q_n\,\colon
n\in\omega\}\subseteq V^P$ such that
$\bfone\forces
``$if $\dot q\in\dot Q\cap
M[G_P]$ then $\dot q_0=\dot q$ and $\dot q_{n+1}\leq\dot q_n$ and
$\dot q_n\in J_n$ and $\delta_n<\sup({\rm range}(\dot q_n))<\delta$.''
Take $\dot r\in V^P$ such that $\bfone\forces``\dot r\in\dot Q$
and ${\rm range}(\dot r)=\bigcup\{{\rm range}(\dot q_n)\,\colon
n\in\omega\}\cup\{\delta\}$.''
Clearly $\dot r$ is as required.

\proclaim Theorem 36 (Shelah). Suppose it is consistent that there is a Mahlo
cardinal. Then it is consistent that GCH holds and
whenever $S\subseteq
S^2_0$ is stationary, then $S$ contains an uncountable
sequentially closed subset.

Proof: The argument given by Shelah in [12, section XI.7]
 can be used, but
with the use of our preservation theorem to simplify
the main lemma there.

As a final comment, let us remark that the proof of
[11, lemma 7]
is incorrect.

\medskip

\medskip

\noindent{\bf References}

\medskip

[1] Baumgartner, J., ``Iterated Forcing'' in {\sl Surveys in Set Theory,}
A.R.D.\  Mathias (ed.), Cambridge University Press, 1979

\medskip

[2] Friedman, H., ``One hundred and two problems in mathematical
logic,'' {\sl Journal of Symbolic Logic,} {\bf 40}, pp.~113--129, 1975

\medskip

[3] Jech, T., {\sl Set Theory,} Academic Press, 1978

\medskip

[4] Jech, T., M. Magidor, W. Mitchell, and K. Prikry, ``Precipitous
ideals,'' {\sl Journal of Symbolic Logic,} vol.~{\bf 45}, pp.~1--8, 1980

\medskip

[5] Jensen, R., and H. Johnsbr\aa ten, ``A new construction of
a non-constructible $\Delta^1_3$ subset of $\omega$,''
{\sl Fundamenta Mathematicae,} vol.~{\bf 81}, pp.~279--290, 1974

\medskip

[6] Kunen, K., {\sl Set Theory, an Introduction to Independence
Proofs,} North-Holland, 1980

\medskip

[7] Kunen, K., and J. Paris, ``Boolean extensions and measurable cardinals,''
{\sl Annals of Mathematical Logic} vol.~{\bf 2}, pp.~359--378, 1971

\medskip

[8] Magidor, M., ``Precipitous ideals and $\Sigma^1_4$ sets,''
{\sl Israel Journal of Mathematics,} vol.~{\bf 35}, pp.~109--134, 1980

\medskip

[9] Namba, K., ``Independence proof of $(\omega,\omega_\alpha)$-distributive
law in complete Boolean algebras,''
{\sl Comment Math.\   Univ.\   St.~Pauli,} vol.~{\bf 19}, pp.~1--12, 1970

\medskip

[10] Schlindwein, C.,
``Consistency of Suslin's hypothesis, a non-special Aronszajn tree,
and GCH,'' {\sl Journal of Symbolic Logic,} vol.~{\bf 59}, pp.~1--29, 1994

\medskip

[11] Schlindwein, C., ``Simplified RCS iterations,''
{\sl Archive for Mathematical Logic,} vol.~{\bf 32}, pp.~341--349, 1993

\medskip

[12] Shelah, S., {\sl Proper Forcing,} Lecture Notes in Mathematics
{\bf 940}, Springer-Verlag, 1982

\vfill\eject

\end{document}